\documentclass[11pt, reqno]{amsart}
\usepackage{amsfonts}
\usepackage{mathrsfs}
\usepackage{amssymb,amsmath,amscd}
\begin{document}
\renewcommand{\a}{\alpha}
\renewcommand{\b}{\beta}
\renewcommand{\d}{\delta}
\newcommand{\D}{\Delta}
\newcommand{\f}{\frac}
\newcommand{\g}{\gamma}
\newcommand{\G}{\Gamma}
\renewcommand{\l}{\lambda}
\renewcommand{\L}{\Lambda}
\newcommand{\be}{\begin{equation}}
\newcommand{\ee}{\end{equation}}
\newcommand{\bea}{\begin{eqnarray}}
\newcommand{\eea}{\end{eqnarray}}
\newcommand{\bna}{\begin{eqnarray*}}
\newcommand{\ena}{\end{eqnarray*}}
\renewcommand{\o}{\omega}
\renewcommand{\O}{\Omega}
\newcommand{\ov}{\over}
\renewcommand{\le}{\left}
\newcommand{\ri}{\right}
\newcommand{\s}{\sigma}
\renewcommand{\th}{\theta}
\newcommand{\ve}{\varepsilon}
\newcommand{\vp}{\varphi}
\renewcommand{\theequation}{\arabic{section}.\arabic{equation}}
\centerline{\bf\Large Asymptotics for cuspidal representations}

\centerline{\bf\Large by functoriality from $GL(2)$}

\bigskip
\centerline{\sc Huixue Lao${}^{1,2}$, Mark McKee${}^2$ {\it and} Yangbo Ye${}^2$}

\centerline{${}^1$ School of Mathematical Sciences, Shandong 
Normal University}

\centerline{Jinan, Shandong 250014, China}

\centerline{${}^2$ Department of Mathematics, The University of Iowa}

\centerline{Iowa City, Iowa 52242-1419, USA}

\centerline{Email: {\tt lhxsdnu@163.com}, 
{\tt huixue-lao@uiowa.edu}, 
{\tt mark-mckee@uiowa.edu},}

\centerline{{\tt yangbo-ye@uiowa.edu}}

\bigskip
\bigskip

\date{}
\openup 1\jot
\bigskip

{\footnotesize {\bf  Abstract} \,\,Let $\pi$ be a unitary 
automorphic cuspidal representation of 
$GL_2(\mathbb{Q}_\mathbb{A})$ with Fourier coefficients 
$\lambda_\pi(n)$. Asymptotic expansions 
of certain sums of $\lambda_\pi(n)$ are proved using 
known functorial liftings from $GL_2$, including 
symmetric powers, isobaric sums, exterior square from 
$GL_4$ and base change. These asymptotic expansions 
are manifestation of the underlying functoriality and 
reflect value distribution of $\lambda_\pi(n)$ on 
integers, squares, cubes and fourth powers.}

\bigskip


\bigskip


\medskip
\section{Introduction}
\setcounter{equation}{0}

Let $\pi=\otimes\pi_p$ be a unitary automorphic irreducible 
cuspidal representation of $GL_2(\mathbb{Q}_\mathbb{A})$ 
with Fourier coefficients $\lambda_\pi(n)$. 
For each finite $p$ for which $\pi_p$ is unramified we can associate to $\pi_p$ Satake parameters 
$\alpha_{\pi} (p), \beta_{\pi} (p) \in \mathbb{C}$ 
such that 
$$ 
\lambda_{\pi} (p) =\alpha_{\pi} (p) +\beta_{\pi} (p),\ \ 
\alpha_{\pi} (p) \beta_{\pi} (p) = 1.
$$
For ramified $\pi_p$,  one can still choose 
$\alpha_{\pi} (p), \beta_{\pi} (p) \in \mathbb{C}$, by 
allowing some of them to be zero, such that the finite-part 
automorphic $L$-function attached to $\pi$ is 
\begin{equation}\label{Lpi}
L(s,\pi) 
=\sum_{n=1}^{\infty} \frac{\lambda_{\pi}(n)}{n^s}
=\prod_{p<\infty}
\big(1-\alpha_\pi(p)p^{-s} \big)^{-1} 
\big(1- \beta_\pi(p)p^{-s} \big)^{-1}
\end{equation}
for Re $s>1$. 

If $\pi$ is associated with a holomorphic cusp form for a 
congruence subgroup of $SL_2(\mathbb Z)$, the Ramanujan 
bound was proved by Deligne \cite{Del} 
\begin{eqnarray*}
|\alpha_\pi (p)|= |\beta_\pi(p)|=1
\end{eqnarray*}
for $p$ with $\pi_p$ being unramified. 
For other cases the best known bound toward the Ramanujan conjecture is
\begin{eqnarray}\label{7over64}
|\alpha_\pi(p)|, |\beta_\pi(p)| \leq p^{\frac{7}{64}},
\end{eqnarray}
which is due to Kim and Sarnak \cite{KimSar}.

Let us review known functoriality from an automorphic cuspidal representation $\pi$ of $GL_2(\mathbb Q_{\mathbb A})$ (cf. 
Kim \cite{Kim2}).

{\it 1.1. The symmetric square lift.} By Gelbart and Jacquet 
\cite{GlbJcq} $Sym^2\pi$ is an automorphic representation 
of $GL_3(\mathbb Q_{\mathbb A})$. It is cuspidal if and only if 
$\pi$ is not monomial. Here $\pi$ is monomial which means that 
there is a nontrivial gr\"ossencharacter $\omega$ so that 
$\pi \cong \pi \otimes \omega$. This in turn means that 
$\pi$ corresponds to a dihedral Galois representation 
Gal$(\bar{\mathbb{Q}}/\mathbb{Q})$. 

{\it 1.2. The symmetric cube.} Kim and Shahidi \cite{KimShh1} 
proved that $Sym^3\pi$ is an automorphic representation 
of $GL_4(\mathbb Q_{\mathbb A})$. 
In addition, $Sym^3\pi$ is cuspidal if and only if 
$\pi$ does not correspond to a dihedral or tetrahedral Galois representation. We recall $\pi$ is tetrahedral means there 
exists a gr\"ossencharacter $\eta \not= 1$ so that 
$Ad^3 \pi \cong Ad^3 \pi \otimes \eta$.  Here 
$Ad^3\pi = Sym^3\pi \otimes \omega_{\pi}^{-1}$, where $\omega_{\pi}$ is the central character of $\pi$. 

{\it 1.3. The symmetric fourth power.} By Kim \cite{Kim1}
$Sym^4\pi$ is an automorphic representation of 
$GL_5(\mathbb Q_{\mathbb A})$. It is either cuspidal or 
unitarily induced from cuspidal representations of 
$GL_2(\mathbb Q_{\mathbb A})$ and 
$GL_3(\mathbb Q_{\mathbb A})$. More precisely, 
$Sym^4\pi$ is cuspidal unless $\pi$ corresponds to 
either a dihedral, tetrahedral, or octahedral Galois 
representation (Kim and Shahidi \cite{KimShh2}). Here 
$\pi$ is octahedral means $Ad^3\pi$ is cuspidal, but there exists a nontrivial quadratic character
$\mu$ so that $Ad^3\pi \cong Ad^3\pi \otimes \mu$.

{\it 1.4. The tensor product of $GL_2\times GL_3$.} 
Let $\pi_2$ and $\pi_3$ be unitary automorphic cuspidal 
representations of $GL_2(\mathbb{Q}_\mathbb{A})$ and 
$GL_3(\mathbb{Q}_\mathbb{A})$, respectively. 
Then, by Kim and Shahidi \cite{KimShh1}, 
$\pi_2 \boxtimes \pi_3$ is an automorphic representation
of $GL_6(\mathbb{Q}_\mathbb{A})$. It is isobaric and is 
cuspidal or irreducibly induced from unitary cuspidal representations. When $\pi_2$ is not dihedral, 
$\pi_2 \boxtimes \pi_3$  is cuspidal unless 
$\pi_3$ is a twist of $Ad(\pi_2)$ by a gr\"ossencharacter 
(Ramakrishnan and Wang \cite{RmkWng}). 

{\it 1.5. The Rankin-Selberg convolution of $GL_2\times GL_2$.} Let $\pi_1$ and $\pi_2$ be unitary automorphic cuspidal 
representations of $GL_2(\mathbb{Q}_\mathbb{A})$. By 
Ramakrishnan \cite{Rmk} 
their isobaric representation $\pi_1\boxtimes\pi_2$ is an 
automorphic representation of $GL_4(\mathbb Q_{\mathbb A})$. 
In case of $\pi_1$ and $\pi_2$ being not associated to characters of a quadratic extension, $\pi_1\boxtimes\pi_2$ 
is cuspidal if and only if $\pi_1$ and $\pi_2$ 
are not twisted equivalent by any gr\"ossencharacter. 

{\it 1.6. Exterior square of $GL_4$.} 
Let $\pi_4$ be a unitary automorphic cuspidal representations 
of $GL_4(\mathbb{Q}_\mathbb{A})$.  By Kim \cite{Kim1} and 
Henniart \cite{Hnn} the exterior square lift 
$\wedge^2(\pi_4)$ is an automorphic representation of 
$GL_6(\mathbb{Q}_\mathbb{A})$. Moreover $\wedge^2(\pi_4)$ is 
of the form 
$Ind\ \sigma_1\otimes\cdots\otimes\sigma_k$  where each $\sigma_i$ is a cuspidal representation of
$GL_{n_i}(\mathbb{Q}_\mathbb{A})$, $n_i >1$. Its cuspidality 
was classified by Asgari and Raghuram \cite{AsgRgh}. 

{\it 1.7. Base change.}  
Let $\pi$ be a unitary automorphic cuspidal representation of $GL_2$ over $\mathbb{Q}$. Suppose that
there is a Galois number field $F$ of prime degree $\ell$ over 
$\mathbb{Q}$ such that
\begin{eqnarray}\label{notpieta}
 \pi \ncong \pi \otimes \eta_{F/\mathbb{Q}},
\end{eqnarray}
where $\eta_{F/\mathbb{Q}}$ is a class character of 
$\mathbb{Q}_{\mathbb{A}}^{\times}$ attached to $F$.
In other words, $\eta_{F/\mathbb{Q}}$ is trivial exactly on 
$\mathbb{Q}^{\times} N_{F/\mathbb{Q}}(F_{\mathbb{A}}^{\times})$. 
Note that $\eta_{F/\mathbb{Q}}$ is a character of order $\ell$. 
By Langlands \cite{Lng} and Arthur and Clozel \cite{ArtClz}, 
\eqref{notpieta} implies that
\begin{eqnarray*}
\pi, \pi \otimes \eta_{F/\mathbb{Q}}, \ldots, 
\eta_{F/\mathbb{Q}}^{\ell-1}
\end{eqnarray*}
lift to an automorphic cuspidal representation $\Pi$ of $GL_2$ 
over $F$. This base change $\Pi$ is stable under the action of 
$Gal (F/\mathbb{Q})$ and given by
\begin{eqnarray}\label{LPi}
L(s, \Pi) 
=\prod_{j=0}^{\ell-1}  
L\big(s, \pi \otimes \eta_{F/\mathbb{Q}}^{j}\big).
\end{eqnarray}

The subject goal of this paper is to establish an asymptotic 
expansion of a certain sum of Fourier coefficient 
$\lambda_\pi(n)$ of the $GL_2$ cuspidal representation 
$\pi$ based on its functorial lifting in each of the cases 
1.1-1.7. These sums and asymptotic expansions can be regarded as 
a manifestation of the underlying functoriality. They also 
provide deep insight into the value distribution of 
$\lambda_\pi(n)$ on positive integers, squares, cubes and 
fourth powers. 

Asymptotics of this kind are important consequences in 
Rankin \cite{Rnk} and Selberg \cite{Slb} as they developed their 
theory. Similar asymptotics were later obtained by 
Lao and Sankaranarayanan \cite{LaoSnk1} \cite{LaoSnk2}, 
Lau, L\"u and Wu \cite{LauLvWu}, 
and Lau and L\"u \cite{LauLv}. 

{\bf Theorem 1.1.} {\it Let $\pi$ be a unitary automorphic 
irreducible cuspidal representation of 
$GL_2(\mathbb{Q}_\mathbb{A})$. Assume that $\pi$ is not 
monomial so that $Sym^2\pi$ is cuspidal. Then the $L$-function 
$$
L_2 (s) 
= 
\sum_{n=1}^{\infty} 
\frac{|\lambda_{\pi}(n^2)|^2}{n^s}
$$ 
converges absolutely for Re $s>1$ and has an analytic 
continuation to Re $s>\frac{1}{2}$ with a simple pole at 
$s=1$. Moreover, 
\begin{equation}\label{n2asym}
\sum_{n \leq x}
|\lambda_{\pi}(n^2)|^2
=
c_{2}x
+O_{\pi}\big(x^{\frac{4}{5}}\big).
\end{equation}
where $c_2 = {\mathop{Res}\limits_{s=1}}\, L_2 (s)$.}

{\bf Theorem 1.2.} {\it Let $\pi$ be a unitary automorphic 
irreducible cuspidal representation of 
$GL_2(\mathbb{Q}_\mathbb{A})$. Assume that $\pi$ does not 
correspond to a dihedral or tetrahedral Galois 
representation so that $Sym^3\pi$ is cuspidal. Then the 
$L$-function 
$$
L_3 (s) 
= 
\sum_{n=1}^{\infty} 
\frac{|\lambda_\pi(n^3)|^2}{n^s}
$$ 
converges absolutely for Re $s>1$, and has an analytic 
continuation to Re $s> \frac{21}{32}$ with a simple pole at 
$s=1$. Moreover, 
$$
\sum_{n \leq x}
|\lambda_\pi(n^3)|^2
=
c_{3}x
+O_{\pi}\big(x^{\frac{15}{17}}\big),
$$
where $c_3 = {\mathop{Res}\limits_{s=1}}\,  L_3 (s)$.}

{\bf Theorem 1.3.} {\it Let $\pi$ be a unitary automorphic 
irreducible cuspidal representation of 
$GL_2(\mathbb{Q}_\mathbb{A})$. Assume that $\pi$ does not 
correspond to a dihedral, tetrahedral, or octahedral 
Galois representation so that $Sym^4\pi$ is cuspidal. 
Then the $L$-function 
$$
L_4 (s) 
= 
\sum_{n=1}^{\infty} 
\frac{|\lambda_\pi(n^4)|^2}{n^s}
$$ 
converges absolutely for Re $s>1$ and has an analytic 
continuation  to Re $s> \frac{7}{8}$ with a simple pole 
at $s=1$. Moreover,
$$
\sum_{n \leq x}
|\lambda_\pi(n^4)|^2
=
c_{4}x
+O_\pi\big(x^{\frac{12}{13}}\big),
$$
where $c_4 = {\mathop{Res}\limits_{s=1}}\,  L_4 (s)$.}

{\bf Theorem 1.4.} {\it Let $\pi_1$ and $\pi_2$ be two 
unitary automorphic cuspidal representations of 
$GL_2(\mathbb{Q}_\mathbb{A})$. Suppose that 
$\pi_1 \boxtimes Sym^2(\pi_2)$ is cuspidal. 
Then the $L$-function 
$$
L_{1,2} (s) 
= 
\sum_{n=1}^{\infty} 
\frac{|\lambda_{\pi_1}(n)|^2
|\lambda_{\pi_2}(n^2)|^2}{n^s}
$$ 
converges absolutely for Re $s>1$ and has an analytic 
continuation to Re $s> \frac{21}{32}$ with a simple pole 
at $s=1$. Moreover,   
\begin{eqnarray}\label{Thm14}
\sum_{n \leq x} 
|\lambda_{\pi_1}(n)|^2
|\lambda_{\pi_2}(n^2)|^2 
=c_{1,2}x
+O_{\pi_1,\pi_2}\big(x^{\frac{35}{37}}\big),
\end{eqnarray}
where $c_{1,2} = {\mathop{Res}\limits_{s=1}}\,  L_{1,2} (s)$.}

{\bf Theorem 1.5.} {\it Let $\pi_1$ and $\pi_2$ be two 
unitary automorphic cuspidal representations of 
$GL_2(\mathbb{Q}_\mathbb{A})$ not associated to characters 
from a quadratic field. Suppose that $\pi_1$ and 
$\pi_2$ are not twisted equivalent by any gr\"ossencharacter 
so that $\pi_1\boxtimes\pi_2$ is cuspidal. Then the $L$-function 
$$
L_{1,1} (s) 
= 
\sum_{n=1}^{\infty} 
\frac{|\lambda_{\pi_1}(n)|^2
|\lambda_{\pi_2}(n)|^2}{n^s}
$$ 
converges absolutely for Re $s>1$ and has an analytic 
continuation to Re $s> \frac12$ with a simple pole 
at $s=1$. Moreover,   
\begin{eqnarray}\label{Thm15}
\sum_{n \leq x} 
|\lambda_{\pi_1}(n)|^2
|\lambda_{\pi_2}(n)|^2 
=c_{1,1}x
+O_{\pi_1,\pi_2}\big(x^{\frac{15}{17}}\big),
\end{eqnarray}
where 
$c_{1,1} = {\mathop{Res}\limits_{s=1}}\,L_{1,1}(s)$.}

{\bf Theorem 1.6.} {\it Let $\pi_1$ and $\pi_2$ be two 
unitary automorphic cuspidal representations of 
$GL_2(\mathbb{Q}_\mathbb{A})$. Suppose that $\pi_1$ and 
$\pi_2$ are not twisted equivalent by any gr\"ossencharacter 
and that $\pi_1\boxtimes\pi_2$ is cuspidal. Then 
\begin{eqnarray}\label{Thm16}
&&
\sum_{abcd \leq x} 
|\lambda_{\pi_1}(a^2)|^2
\lambda_{\pi_1}(b^2)\overline{\lambda}_{\pi_2}(b^2)
\lambda_{\pi_2}(c^2)\overline{\lambda}_{\pi_1}(c^2)
|\lambda_{\pi_2}(d^2)|^2 
\nonumber
\\
&=&
x(c_1\log x+c_0)
+O_{\pi_1,\pi_2}(x^{\frac{35}{37}}\log x)
\end{eqnarray}
for some constants $c_1$ and $c_0$.}

{\bf Theorem 1.7.} {\it Let $\pi$ be a unitary automorphic 
cuspidal representation of $GL_2(\mathbb{Q}_{\mathbb{A}})$, 
and $\chi$ a nontrivial Dirichlet character modulo prime 
$\ell$. Then
\begin{equation}\label{Thm17}
\sum_{\substack{a_{ij} \geq 1 \\ 
\prod_{i=0}^{\ell-1} \prod_{j=0}^{\ell-1} a_{ij} \leq x}}
\prod_{i=0}^{\ell-1} \prod_{j=0}^{\ell-1} 
|\lambda_{\pi} (a_{ij})|^2 \chi^j (a_{ij}) 
= 
x P_{\ell-1} (\log x) 
+O_{\pi,\chi} \Big(x^{\frac{4 \ell^2 -1}{4 \ell^2 +1}} 
\log^{\ell-1} x \Big),
\end{equation}
where $P_{\ell-1}$ is a polynomial of degree $\ell-1$.}

We remark that Theorems 1.1-1.3 
are a generalization of Theorem 1.1 in Lao and Sankaranarayanan 
\cite{LaoSnk1} \cite{LaoSnk2} 
where the same asymptotics (with the same error terms) 
were proved assuming $\pi$ corresponds to a self-dual 
holomorphic cusp form of even integral weight 
on the full modular group $SL_2(\mathbb Z)$. The proof of 
Theorem 1.2 in \cite{LaoSnk2}, however, should not have used 
a subconvexity bound for $L(s,Sym^2f)$ proved by Li \cite{XLi}. 
Using a subconvexity bound for $\zeta(s)$ which was recently 
improved by Bourgain \cite{Brg} we can correct the proof and 
obtain a better result.

{\bf Theorem 1.8.} {\it Let $f$ be a self-dual 
holomorphic cusp form 
for a congruence subgroup of $SL_2(\mathbb Z)$. Then}
\begin{equation}\label{n2asymholo}
\sum_{n \leq x}
\lambda^2_{f}(n^2)
=
c_{2}x
+O_f\big(x^{\frac{566}{737}}\big).
\end{equation}

In the next section we will outline the main techniques which 
we will use. 

\section{Some properties of $L$-functions}
\setcounter{equation}{0}

We need a classical lemma by Landau \cite{Lnd} (cf. Ram Murty \cite{Mrt}). 

{\bf Lemma 2.1.} {\it Let $a_n\geq 0$ and set
$$f(s)=\sum_{n=1}^{\infty} \frac{a_n}{n^s}.$$
Suppose $f(s)$ is convergent
in some half-plane and that it has analytic continuation, except
for a pole at $s=1$ of order $k$, to the entire complex plane
and it satisfies a functional equation
$$
c^s\Delta(s)f(s)=c^{1-s}\Delta(1-s)f(1-s),
$$
where $c$ is a certain positive constant and 
$$
\Delta(s)=\prod_{i=1}^N
\Gamma(\alpha_i s+\beta_i) (\alpha_i > 0).
$$
Then
$$
\sum_{n \leq x}a_n=x P_{k-1}(\log x)+O\Big(x^{\frac{2A-1}{2A+1}}\log^{k-1}x\Big),
$$
where $A=\sum_{i=1}^N\alpha_i$, and $P_{k-1} (y)$ is a polynomial in $y$ of degree $k-1$.}

Now we summarize some basic properties of automorphic 
$L$-functions. Let $\pi$ be a unitary automorphic cuspidal 
representation of $GL_2(\mathbb Q_{\mathbb A})$,  with  
Fourier coefficients $\lambda_{\pi}(n)$ normalized by 
$\lambda_\pi (1)=1$. For simplification, we write 
$\alpha_{\pi}(p)=\alpha_p$ and $\beta_{\pi} (p)=\beta_p$. 
By \eqref{Lpi}, we have
\begin{equation}\label{lambdapk} 
\lambda_\pi(p^k)
= \sum_{i=0}^k 
\alpha_p^{k-i} \beta_p^i. 
\end{equation}

Let $S$ be a finite set of primes containing all $p$ at which 
$\pi_p$ is ramified. Then the $m^{\text{th}}$-symmetric power 
$L$-function attached to $\pi$ 
$$
L(s, Sym^m \pi)
=
\sum_{n=1}^{\infty}
\frac{\lambda_{{Sym}^m\pi} (n)}{n^s}
$$
has an Euler product outside $S$
\begin{eqnarray*}
L^S(s, Sym^m \pi)
=
\prod_{p\notin S}\prod_{i=0}^m
(1-\alpha_p^{m-i}\beta_p^ip^{-s})^{-1}
\end{eqnarray*}
for Re s $>1$. 

Similarly the Rankin-Selberg $L$-function of $Sym^m\pi$ times 
its contragredient 
$$
L(s, Sym^m\pi \times Sym^m\widetilde\pi)
=
\sum_{n=1}^{\infty}
\frac{\lambda_{Sym^m\pi\times Sym^m\widetilde\pi}(n)}{n^s}
$$
has an Euler product outside $S$
$$
L^S(s, Sym^m\pi \times Sym^m\widetilde\pi)
=
\prod_{p\notin S}\prod_{i=0}^m\prod_{j=0}^m
(1-\alpha_p^{m-i}\beta_p^i
\bar\alpha_p^{m-j}\bar\beta_p^j p^{-s})^{-1}
$$
for Re $s>1$. For $p\notin S$ its parameters are
\begin{eqnarray*}
  |\alpha_p|^{2m},
\ |\alpha_p|^{2m-2}{\alpha_p}{\bar\beta_p},
\ |\alpha_p|^{2m-2}{\bar\alpha_p}{\beta_p},
\ |\alpha_p|^{2m-4}{\alpha^2_p}{\bar\beta^2_p},
\ |\alpha_p|^{2m-4},
\ |\alpha_p|^{2m-4}{\bar\alpha^2_p}{\beta^2_p},
\ \ldots, 
\\
  |\beta_p|^{2m-4}{\beta^2_p}{\bar\alpha^2_p},
\ |\beta_p|^{2m-4},
\ |\beta_p|^{2m-4}{\bar\beta^2_p}{\alpha^2_p},
\ |\beta_p|^{2m-2}{\beta_p}{\bar\alpha_p},
\ |\beta_p|^{2m-2}{\bar\beta_p}{\alpha_p},
\ |\beta_p|^{2m}.
\end{eqnarray*}
Consequently 
\begin{eqnarray*}\label{RSsymmExp}
&&
L^S(s,Sym^m \pi \times Sym^m \widetilde\pi)
\nonumber
\\
&=& 
\prod_{p\notin S} 
\Big(1-\frac{|\alpha_p|^{2m}}{ p^{s}}\Big)^{-1} 
\Big(1-\frac{ |\alpha_p|^{2m-2}{\alpha_p}{\bar\beta_p}}
{ p^{s}}\Big)^{-1} 
\Big(1-\frac{ |\alpha_p|^{2m-2}{\bar\alpha_p}{\beta_p}}
{ p^{s}}\Big)^{-1} 
\nonumber
\\
&&\times
\Big(1-\frac{ |\alpha_p|^{2m-4}{\alpha_p^2}
{\bar\beta_p^2}}{ p^{s}}\Big)^{-1} 
\Big(1-\frac{ |\alpha_p|^{2m-4}}{ p^{s}}\Big)^{-1} 
\Big(1-\frac{ |\alpha_p|^{2m-4}{\bar\alpha_p^2}
{\beta_p^2}}{ p^{s}}\Big)^{-1} 
\cdots.
\end{eqnarray*}
Then we have for $p\notin S$ 
\begin{eqnarray}\label{lambdaRSsymm}
&&
\lambda_{Sym^m  \pi \times Sym^m  \widetilde\pi} (p^{j})
\nonumber
\\
&=&
|\alpha_p|^{2mj} + 
|\alpha_p|^{2mj-2} {\alpha_p}{\bar\beta_p} +
|\alpha_p|^{2mj-2} {\bar\alpha_p}{\beta_p} 
\nonumber
\\
&&+
2|\alpha_p|^{2mj-4} {\alpha^2_p}{\bar\beta^2_p} +
2|\alpha_p|^{2mj-4} +
2|\alpha_p|^{2mj-4} {\bar\alpha^2_p}{\beta^2_p} 
+\cdots
\end{eqnarray}
for $j\geq1$.

We note that for $m=2,3,4$, the Rankin-Selberg $L$-function 
$L(s, Sym^m\pi\times Sym^m\widetilde\pi)$ has an analytic continuation to $\mathbb{C}$ with a simple pole at $s=1$, 
when $Sym^m \pi$ is cuspidal. Using the convexity bound for 
$L(s, Sym^4\pi\times Sym^4\widetilde\pi)$ Li and Young 
\cite{XLiYng} proved the following lemma.

{\bf Lemma 2.2.} {\it Let $\pi$ be a unitary automorphic 
cuspidal representation of $GL_2(\mathbb Q_{\mathbb A})$ such 
that $Sym^4\pi$ is cuspidal. Then 
$$
\prod_p
\Big(
1+\frac{|\alpha_p|^8+|\beta_p|^8}{p^\sigma}
\Big)
$$
converges for $\sigma>1$.}

\section{ Proof of Theorems 1.1-1.3.}
\setcounter{equation}{0}

Define a Dirichlet series
\begin{equation}\label{Lm}
L_m(s)
=
\sum_{n=1}^{\infty}\frac{|\lambda_\pi(n^m)|^2}{n^s}
=
\prod_{p} 
\Big(
\sum_{j=0}^{\infty} 
\frac{|\lambda_\pi(p^{mj})|^2}{p^{js}} 
\Big)
\end{equation}
for Re $s$ sufficiently large. 
Since $L_m (s)$ is not known to be automorphic, we want to 
factor a Rankin-Selberg $L$-function out of $L_m (s)$.

{\bf Lemma 3.1} {\it Let $\pi$ be a unitary automorphic 
cuspidal representation of $GL_2(\mathbb Q_{\mathbb A})$ such 
that $Sym^m\pi$ is cuspidal.
Then we have for $m=2,3,4$,
\begin{eqnarray*}
L_m(s)=L(s, Sym^m\pi \times Sym^m \widetilde\pi)U_m(s),
\end{eqnarray*}
where $U_m(s)$ as an Euler product is absolutely convergent 
and hence nonzero for Re $s>\frac{1}{2}$ for $m=2$, 
Re $s>\frac{23}{32}$ for $m=3$, Re $s>\frac{7}{8}$ for 
$m=4$.}

{\it Proof.}
First we factor $L_m(s)$ as in \eqref{Lm}:
\begin{eqnarray}\label{LmUm}
L_m(s)
&=&
\prod_p 
\Big(
\sum_{j=0}^{\infty} 
\frac{\lambda_{Sym^m\pi\times Sym^m\widetilde\pi}(p^j)}
{p^{js}}
-
\sum_{j=1}^{\infty} 
\frac{P_j (\alpha_p, \beta_p)}{p^{js}}
\Big)
\nonumber
\\
&=& 
\prod_p 
\Big(
\sum_{j=0}^{\infty} 
\frac{\lambda_{Sym^m\pi\times Sym^m\widetilde\pi}(p^j)}
{p^{js}}
\Big) 
\prod_{p} 
\Big( 
1- \frac{\sum_{j=1}^{\infty} P_j (\alpha_p, \beta_p) p^{-js}}{\sum_{j=0}^{\infty}
\lambda_{Sym^m\pi\times Sym^m\widetilde\pi}(p^j) p^{-js}}   \Big)
\nonumber
\\ 
&=& 
L(s,Sym^m  \pi \times Sym^m \widetilde\pi ) 
U_m(s)
\end{eqnarray}
by denoting the last product in \eqref{LmUm} by $U_m(s)$. 
Here 
$$
P_j (\alpha_p, \beta_p)
= 
\lambda_{Sym^m\pi\times Sym^m\widetilde\pi}(p^j)
-
|\lambda_\pi(p^{mj})|^2
$$
for $j \geq 1$.

By \eqref{lambdapk}, we have 
\begin{eqnarray}\label{lambdapi2}
|\lambda_\pi ( p^{mj})|^2
&=&
\Big|\sum_{i=0}^{mj} \alpha_p^{mj-i} \beta_p^i \Big|^2 
\nonumber
\\
&=&
|\alpha_p|^{2mj} + 
|\alpha_p|^{2mj-2} {\alpha_p}{\bar\beta_p} +
|\alpha_p|^{2mj-2} {\bar\alpha_p}{\beta_p} 
\nonumber
\\
&&+
|\alpha_p|^{2mj-4} {\alpha^2_p}{\bar\beta^2_p} +
|\alpha_p|^{2mj-4} +
|\alpha_p|^{2mj-4} {\bar\alpha^2_p}{\beta^2_p} 
+\cdots
\end{eqnarray}
for $j\geq1$. Comparing \eqref{lambdaRSsymm} and 
\eqref{lambdapi2} we get for $p\notin S$ 
\begin{eqnarray*}
P_1 (\alpha_p, \beta_p)
&=&
0, 
\\
P_j(\alpha_p, \beta_p)
&=& 
|\alpha_p|^{2mj-4} {\alpha^2_p}{\bar\beta^2_p} +
|\alpha_p|^{2mj-4} +
|\alpha_p|^{2mj-4} {\bar\alpha^2_p}{\beta^2_p} 
+\cdots
\\
&&+
|\beta_p|^{2mj-4} {\beta^2_p}{\bar\alpha^2_p} +
|\beta_p|^{2mj-4} +
|\beta_p|^{2mj-4} {\bar\beta^2_p}{\alpha^2_p}.
\end{eqnarray*}
Now let us turn to the series in the denominator in 
\eqref{LmUm} for $p\notin S$. By \eqref{RSsymmExp}, we get
\begin{eqnarray*}
&& 
\Big(\sum_{j=0}^{\infty} 
\frac{\lambda_{Sym^m\pi\times Sym^m\widetilde\pi}(p^j)}
{p^{js}}\Big)^{-1}
\\
&=& 
\Big(1-\frac{|\alpha_p|^{2m}}{ p^{s}}\Big) 
\Big(1-\frac{ |\alpha_p|^{2m-2}{\alpha_p}{\bar\beta_p}}
{ p^{s}}\Big) 
\Big(1-\frac{ |\alpha_p|^{2m-2}{\bar\alpha_p}{\beta_p}}
{ p^{s}}\Big) 
\nonumber
\\
&&\times
\Big(1-\frac{ |\alpha_p|^{2m-4}{\alpha_p^2}
{\bar\beta_p^2}}{ p^{s}}\Big) 
\Big(1-\frac{ |\alpha_p|^{2m-4}}{ p^{s}}\Big) 
\Big(1-\frac{ |\alpha_p|^{2m-4}{\bar\alpha_p^2}
{\beta_p^2}}{ p^{s}}\Big) 
\cdots.
\end{eqnarray*}
Thus for $p\notin S$
\begin{eqnarray}\label{Umterm}
&&
 - \frac{
\sum_{j=2}^{\infty} P_j (\alpha_p, \beta_p) p^{-js}
}
  {
\sum_{j=0}^{\infty} 
\lambda_{Sym^m\pi\times Sym^m\widetilde\pi}(p^j) p^{-js}
} 
\nonumber
\\
&\ll&  
\frac{|\alpha_p|^{4m-4} +|\beta_p|^{4m-4}}{p^{2 \sigma}}+
  \frac{|\alpha_p|^{6m-4} +|\beta_p|^{6m-4}}{p^{3 \sigma}}+ \cdots.
\end{eqnarray}

First let us consider the case of $m=2$. Then for 
$p\notin S$ \eqref{Umterm} 
becomes
$$ 
\ll  
\frac{|\alpha_p|^{4} +|\beta_p|^{4}}{p^{2 \sigma}}+
  \frac{|\alpha_p|^{8} +|\beta_p|^{8}}{p^{3 \sigma}}+ \cdots
$$
which converges when $\sigma > \frac{7}{16}$. Thus $U_2(s)$ in 
\eqref{LmUm} is dominated by
\begin{eqnarray}\label{dominateU2}
\prod_{p\notin S} 
\Big(1+\frac{|\alpha_p|^{4} +|\beta_p|^{4}}{p^{2 \sigma}} \Big).
\end{eqnarray}
By Cauchy's inequality, \eqref{dominateU2} is dominated 
in turn by
$$
\prod_{p\notin S} 
\Big(1+\frac{|\alpha_p|^{8} +|\beta_p|^{8})}
{p^{4 \sigma -1-\varepsilon }} \Big)^{\frac12} 
\prod_{p\notin S} \Big( 1+ \frac{1}{p^{1+\varepsilon}}\Big)^{\frac12} 
$$
which converges when $4 \sigma -1-\varepsilon >1$, i.e., 
when $\sigma > \frac12 + \varepsilon,$  by Lemma 2.2.
Hence $U_2(s)$ is absolutely convergent and hence is nonzero for  $\sigma > \frac{1}{2}$.

Then let us turn to the case of $m=3$. Then for $p\notin S$ 
\eqref{Umterm} 
becomes
$$ 
\ll  
\frac{|\alpha_p|^{8} +|\beta_p|^{8}}{p^{2 \sigma}}+
  \frac{|\alpha_p|^{14} +|\beta_p|^{14}}{p^{3 \sigma}}+ \cdots
$$
  which converges when $\sigma >  \frac{21}{32}$. Thus $U_3(s)$ in \eqref{LmUm} is dominated by
$$
\prod_{p\notin S} 
\Big(1+\frac{|\alpha_p|^{8} +|\beta_p|^{8})}{p^{2 \sigma}} \Big)
$$
which converges absolutely for $\sigma > \frac12$ by Lemma 2.2.
Therefore, $U_3(s)$ converges absolutely for $\sigma > \frac{21}{32}.$

Now let $m=4$. Then for $p\notin S$ \eqref{Umterm} is
$$ 
\ll  
\frac{|\alpha_p|^{12} +|\beta_p|^{12}}{p^{2 \sigma}}+
  \frac{|\alpha_p|^{20} +|\beta_p|^{20}}{p^{3 \sigma}}+ \cdots
$$
 which converges for $\sigma >  \frac{7}{8}$. Thus $U_4(s)$ in 
\eqref{LmUm} is dominated by
$$
\prod_{p\notin S} 
\Big(1+\frac{|\alpha_p|^{12} +|\beta_p|^{12})}{p^{2 \sigma}} \Big)
$$
which in turn is dominated by
\begin{eqnarray}\label{dominateU4}
\prod_{p\notin S} \Big(1+\frac{|\alpha|^{8} +|\beta|^{8}}{p^{ 2 \sigma - \frac{7}{16}}}  \Big)
\end{eqnarray}
by \eqref{7over64}. 
By Lemma 2.2 \eqref{dominateU4} converges absolutely for 
$2 \sigma - \frac{7}{16}>1,$ i.e., for $\sigma > \frac{23}{32}$. 
Consequently, $U_4(s)$ converges absolutely for Re $s> \frac78$. \hfill 
$\square$

The $L$-function 
$L(s, Sym^m \pi \times Sym^m \widetilde\pi)$ satisfies 
the conditions of Landau's Lemma 2.1 with $k=1$ and 
$2A=(m+1)^2$. Thus we have
\begin{equation}\label{Landaum}
\sum_{n \leq x} 
\lambda_{Sym^m \pi \times Sym^m \widetilde\pi}(n) 
=d_{m}x+O\big(x^{1-\frac{2}{(m+1)^2 +1}}\big),
\end{equation}
where 
$$
d_m={\mathop{Res}\limits_{s=1}}\,  
L(s, Sym^m \pi \times Sym^m \widetilde\pi).
$$
Expand 
$$
U_m(s) 
= \sum_{\ell=1}^{\infty} \frac{u_m(\ell)}{\ell^s}
$$
which converges absolutely for Re $s>\frac12$ for $m=2$, 
Re $s>\frac{23}{32}$ for $m=3$, and Re $s>\frac78$ for $m=4$ 
by Lemma 3.1. By \eqref{LmUm}
\begin{equation}\label{conv}
|\lambda_\pi(n^m)|^2
=
\sum_{n=k\ell}
\lambda_{Sym^m \pi \times Sym^m \widetilde\pi}(k)u_m(\ell).
\end{equation}

For $m=2$, \eqref{conv} becomes
\begin{eqnarray}\label{2n2}
\sum_{n \leq x}
|\lambda_\pi(n^2)|^2
&=&
\sum_{\ell\leq x}
u_2(\ell)
\sum_{k \leq x/\ell}
\lambda_{Sym^2 \pi \times Sym^2 \widetilde\pi}(k)
\nonumber
\\
&=&
\sum_{\ell \leq x}
u_2(\ell) 
\Big\{d_2  \frac{x}{\ell}
+O\Big(\Big(\frac{x}{\ell}\Big)^{\frac{4}{5}}\Big) \Big\}
\nonumber
\\ 
&=& 
x d_2 
\sum_{\ell \leq x} 
\frac{u_2(\ell)}{\ell}
+O\Big( x^{\frac{4}{5}}  
\sum_{\ell \leq x} \frac{|u_2(\ell)|}{\ell^{\frac{4}{5}}} \Big)
\end{eqnarray}
by \eqref{Landaum}. Using 
\begin{eqnarray*}
\sum_{\ell > x}  \frac{u_2(\ell)}{\ell}   
\ll x^{-\frac{1}{2}+\epsilon}, 
\quad
\sum_{\ell \leq x} 
\frac{u_2(\ell)}{\ell^{\frac{4}{5}}} 
\ll 1,
\end{eqnarray*}
we get
\begin{eqnarray}\label{u2}
\sum_{\ell \leq x}  
\frac{u_2(\ell)}{\ell} 
=U_2 (1) 
+ O\big(  x^{-\frac{1}{2}+\epsilon} \big).
\end{eqnarray}
By \eqref{2n2} and \eqref{u2} we conclude
\begin{eqnarray*}
\sum_{n \leq x}\lambda^2(n^2) = c_2 x + O\big(x^{\frac{4}{5}}\big),
\end{eqnarray*}
where $c_2=d_2 U_2 (1)={\mathop{Res}\limits_{s=1}}\,  L_2(s)$ or by a 
Tauberian argument.

The cases of $m=3,4$ are similar.
\hfill $\square$

\section{ Proof of Theorems 1.4-1.7.}
\setcounter{equation}{0}

{\it Proof of Theorem 1.4.} Let $\pi_1$ and $\pi_2$ be 
unitary automorphic cuspidal representations of 
$GL_2(\mathbb{Q}_\mathbb{A})$. Suppose that $\pi_1$ and 
$\pi_2$ are not twisted equivalent by a gr\"ossencharacter 
so that $\pi_1 \boxtimes Sym^2(\pi_2)$ is cuspidal. Denote 
by $ \alpha_p$ and $ \beta_p$ the Langlands parameters of 
$\pi_1$ at $p$, and by $\gamma_p$ 
and $\delta_p$ the Langlands parameters of $\pi_2$. 
Suppose without loss of generality that $|\alpha_p|\geq1$ 
and $\gamma_p|\geq1$.
Let $T$ be a finite set of primes containing all $p$ at which 
either $\pi_{1p}$ or $\pi_{2p}$ is ramified. 
Then at $p\notin T$
$Sym^2\pi_2$ has parameters $\gamma_p^2$, $1$ and $\delta_p^2$, 
while $\pi_1\boxtimes Sym^2\pi_2$ has 
$\alpha_p\gamma_p^2$, $\alpha_p$, $\alpha_p\delta_p^2$, 
$\beta_p\gamma_p^2$, $\beta_p$ and $\beta_p\delta_p^2$. Thus 
for $p\notin T$ 
\begin{equation}\label{12p}
\lambda_{ (\pi_1\boxtimes Sym^2\pi_2) 
\times(\widetilde\pi_1\boxtimes Sym^2\widetilde\pi_2) }(p)
=
|\alpha_p+\beta_p|^2
|\gamma_p^2+1+\delta_p^2|^2
=
|\lambda_{\pi_1}(p)|^2
|\lambda_{\pi_2}(p^2)|^2
\end{equation}
and both 
$\lambda_{ (\pi_1\boxtimes Sym^2\pi_2) 
\times(\widetilde\pi_1\boxtimes Sym^2\widetilde\pi_2) }(p^j)$ 
and 
$|\lambda_{\pi_1}(p^j)|^2|\lambda_{\pi_2}(p^{2j})|^2$ 
for $j\geq2$ are equal to
\begin{eqnarray}\label{both}
&&
|\alpha_p|^{2j}|\gamma_p|^{4j}
+
|\alpha_p|^{2j-4}|\gamma_p|^{4j}
(\alpha_p^2+\bar\alpha_p^2)
+
|\alpha_p|^{2j}|\gamma_p|^{4j-4}
(\gamma_p^2+\bar\gamma_p^2)
\nonumber
\\
&+&
O\Big(
|\alpha_p|^{2j}|\gamma_p|^{4j-4}
+|\alpha_p|^{2j-2}|\gamma_p|^{4j-2}
+|\alpha_p|^{2j-4}|\gamma_p|^{4j}
\Big).
\end{eqnarray}

Define 
\begin{equation}\label{12pj}
Q_j(\alpha_p,\beta_p,\gamma_p,\delta_p)
=
\lambda_{ (\pi_1\boxtimes Sym^2\pi_2) 
\times(\widetilde\pi_1\boxtimes Sym^2\widetilde\pi_2) }(p^j)
-
|\lambda_{\pi_1}( p^j)|^2 
|\lambda_{\pi_2}( p^{2j})|^2
\end{equation}
for $j\geq1$. Then by \eqref{12p} and \eqref{both} we have at 
$p\notin T$ that 
\begin{eqnarray}\label{Qj}
Q_1(\alpha_p,\beta_p,\gamma_p,\delta_p)
&=&
0
\nonumber
\\
Q_j(\alpha_p,\beta_p,\gamma_p,\delta_p)
&\ll&
|\alpha_p|^{2j}|\gamma_p|^{4j-4}
+|\alpha_p|^{2j-2}|\gamma_p|^{4j-2}
+|\alpha_p|^{2j-4}|\gamma_p|^{4j}
\end{eqnarray}
for $j\geq2$. 
Consequently by \eqref{12pj}, the $L$-function
\begin{eqnarray*}
L_{1,2} (s) 
= 
\sum_{n=1}^{\infty} 
\frac{|\lambda_{\pi_1}(n)|^2
|\lambda_{\pi_2}(n^2)|^2}{n^s}
=
\prod_p 
\sum_{j=0}^{\infty} 
\frac{|\lambda_{\pi_1}(p^j)|^2
|\lambda_{\pi_2}(p^{2j})|^2}{p^{js}},
\end{eqnarray*}
for Re $s$ sufficiently large, can be factored as
\begin{eqnarray}\label{L12factor}
L_{1,2} (s) 
&=& 
\prod_p 
\Big(
\sum_{j=0}^{\infty} 
\frac{
\lambda_{ (\pi_1\boxtimes Sym^2\pi_2) 
\times(\widetilde\pi_1\boxtimes Sym^2\widetilde\pi_2) }(p^j)
}
{p^{js}} 
- 
\sum_{j=1}^{\infty} 
\frac{
Q_j(\alpha_p,\beta_p,\gamma_p,\delta_p)
}{p^{js}} 
\Big) 
\nonumber
\\
&=& 
\prod_p 
\Big(
\sum_{j=0}^{\infty} 
\frac{
\lambda_{ (\pi_1\boxtimes Sym^2\pi_2) 
\times(\widetilde\pi_1\boxtimes Sym^2\widetilde\pi_2) }(p^j)
}
{p^{js}} 
\Big)
\nonumber
\\
&&\times
\prod_{p} 
\Big( 1- 
\frac{ \sum_{j=1}^{\infty} 
Q_j (\alpha_p, \beta_p,\gamma_p,\delta_p) p^{-js} }
{\sum_{j=0}^{\infty} 
\lambda_{ (\pi_1\boxtimes Sym^2\pi_2) 
\times(\widetilde\pi_1\boxtimes Sym^2\widetilde\pi_2) }(p^j)
 p^{-js}}   
\Big)
\nonumber
\\
&=&  
L(s,(\pi_1\boxtimes Sym^2\pi_2) 
\times(\widetilde\pi_1\boxtimes Sym^2\widetilde\pi_2)) 
V_{1,2}(s)
\end{eqnarray}
by denoting the last product in \eqref{L12factor} by 
$V_{1,2}(s)$. Note that by \eqref{Qj} the sum of 
$Q_j$ in \eqref{L12factor} is actually taken over $j\geq2$ 
for $p\notin T$. 

For $p\notin T$ by \eqref{Qj} the series 
\begin{equation}\label{dominate}
\sum_{j=2}^{\infty} 
\frac{Q_j (\alpha_p, \beta_p,\gamma_p,\delta_p)}{ p^{js}}
\end{equation}
is dominated by
$$
\sum_{j=2}^{\infty}  
\frac{
|\alpha_p|^{2j}|\gamma_p|^{4j-4}
+|\alpha_p|^{2j-2}|\gamma_p|^{4j-2}
+|\alpha_p|^{2j-4}|\gamma_p|^{4j}
}
{p^{j \sigma}}
$$
which in turn is dominated by
$$
\sum_{j=2}^{\infty}  
\frac{1}{p^{j \sigma-(6j-4) \frac{7}{64}}}
$$
by \eqref{7over64}. Thus \eqref{dominate} converges 
absolutely when $\sigma >  \frac{21}{32}$, and its sum is then
\begin{eqnarray*}
\ll  
\frac{
|\alpha_p|^4|\gamma_p|^4
+|\alpha_p|^2|\gamma_p|^6
+|\gamma_p|^8
}
{p^{2 \sigma}}.
\end{eqnarray*}
 Therefore when $\sigma > \frac{21}{32}$,   $V_{1,2}(s)$ is dominated by
\begin{eqnarray}\label{dominateprod}
\prod_{p\notin T} 
\Big(1+\frac{
|\alpha_p|^4|\gamma_p|^4
+|\alpha_p|^2|\gamma_p|^6
+|\gamma_p|^8
}{p^{2 \sigma}} \Big).
\end{eqnarray}
By Lemma 2.2 and Cauchy's inequality, \eqref{dominateprod} 
converges absolutely for $\sigma > \frac12.$ Thus $V_{1,2} (s)$ 
converges absolutely for Re $s> \frac {21}{32}$. 

Applying Landau's Lemma 2.1 to 
$L(s, (\pi_1\boxtimes Sym^2\pi_2) 
\times(\widetilde\pi_1\boxtimes Sym^2\widetilde\pi_2) )$, 
we have
$$
\sum_{n \leq x} 
\lambda_{
(\pi_1\boxtimes Sym^2\pi_2) 
\times(\widetilde\pi_1\boxtimes Sym^2\widetilde\pi_2)}(n) 
= d_{1,2} x + O \big( x^{\frac{35}{37}}\big).
$$
By the Dirichlet convolution on the right hand side of 
\eqref{L12factor} we have 
\begin{eqnarray*}
\sum_{n \leq x}
|\lambda_{\pi_1}(n)|^2
|\lambda_{\pi_2}(n^2)|^2
&=&
\sum_{\ell\leq x}
u_{1,2}(\ell)
\sum_{k \leq x/\ell}
\lambda_{
(\pi_1\boxtimes Sym^2\pi_2) 
\times(\widetilde\pi_1\boxtimes Sym^2\widetilde\pi_2)}(k) 
\\ 
&=&
\sum_{\ell \leq x}
u_{1,2}(\ell) 
\Big\{d_{1,2}  
\frac{x}{\ell}
+O\Big(\Big(\frac{x}{\ell}\Big)^{\frac{35}{37}}\Big) \Big\},
\end{eqnarray*}
which leads to \eqref{Thm14}.
\hfill $\square$

{\it Proof of Theorem 1.5} is similar to that of Theorem 1.4 
with the same notation. We have for $p\notin T$ 
\begin{eqnarray}\label{11p}
\lambda_{ (\pi_1\boxtimes \pi_2) 
\times(\widetilde\pi_1\boxtimes \widetilde\pi_2) }(p)
&=&
|\alpha_p+\beta_p|^2
|\gamma_p+\delta_p|^2
=
|\lambda_{\pi_1}(p)|^2
|\lambda_{\pi_2}(p)|^2
\\
R_j(\alpha_p,\beta_p,\gamma_p,\delta_p)
&=&
\lambda_{ (\pi_1\boxtimes\pi_2) 
\times(\widetilde\pi_1\boxtimes\widetilde\pi_2) }(p^j)
-
|\lambda_{\pi_1}( p^j)|^2 
|\lambda_{\pi_2}( p^{j})|^2
\nonumber
\\\label{11Rj}
&\ll&
|\alpha_p|^{2j}|\gamma_p|^{2j-4}
+|\alpha_p|^{2j-2}|\gamma_p|^{2j-2}
+|\alpha_p|^{2j-4}|\gamma_p|^{2j}
\end{eqnarray}
for $j\geq2$. 
Consequently by \eqref{11p} and \eqref{11Rj}, the $L$-function
\begin{equation}\label{L11factor}
L_{1,1} (s) 
= 
\sum_{n=1}^{\infty} 
\frac{|\lambda_{\pi_1}(n)|^2
|\lambda_{\pi_2}(n)|^2}{n^s}
=
L(s,(\pi_1\boxtimes\pi_2) 
\times(\widetilde\pi_1\boxtimes\widetilde\pi_2)) 
V_{1,1}(s)
\end{equation}
with $V_{1,1} (s)$ 
convergent absolutely for Re $s> \frac12$. 

Applying Landau's Lemma 2.1 to 
$L(s, (\pi_1\boxtimes\pi_2) 
\times(\widetilde\pi_1\boxtimes\widetilde\pi_2) )$, 
we have
$$
\sum_{n \leq x} 
\lambda_{
(\pi_1\boxtimes\pi_2) 
\times(\widetilde\pi_1\boxtimes\widetilde\pi_2)}(n) 
= d_{1,1} x + O \big( x^{\frac{15}{17}}\big).
$$
By the Dirichlet convolution on the right hand side of 
\eqref{L11factor} we have \eqref{Thm15}.
\hfill$\square$

{\it Proof of Theorem 1.6.} With $\pi_1$ and $\pi_2$ as above 
so that $\pi_1\boxtimes\pi_2$ is cuspidal, its exterior 
square representation $\wedge^2 (\pi_1 \boxtimes \pi_2)$ has Langlands parameters $\alpha_p^2$, $\beta_p^2$, 
$\gamma_p^2$, $\delta_p^2$ and $1\times2$. Thus
$$
\wedge^2 (\pi_1 \boxtimes \pi_2)  
\cong 
Ind (Sym^2 \pi_1 \otimes Sym^2 \pi_2)
$$
and
\begin{eqnarray}\label{exteriorL}
&&
L(s, \wedge^2 (\pi_1\boxtimes \pi_2) \times 
\wedge^2 (\widetilde{\pi}_1 \boxtimes \widetilde{\pi}_2))
\nonumber
\\
&=&
L(s,Sym^2 \pi_1 \times Sym^2 \widetilde{\pi}_1) 
L(s,Sym^2 \pi_1 \times Sym^2 \widetilde{\pi}_2)
\nonumber
\\
&&\times
L(s,Sym^2 \pi_2 \times Sym^2 \widetilde{\pi}_1) 
L(s,Sym^2 \pi_2 \times Sym^2 \widetilde{\pi}_2).
\end{eqnarray}
Since $\pi_1$ and $\pi_2$ are not twisted equivalent by 
any gr\"ossencharacter, \eqref{exteriorL} has a double pole at $s=1$ and is otherwise analytic.
The $L$-function we are interested in is
\begin{eqnarray}\label{4L}
L_{\wedge^2}(s)
&=&
L_{\pi_1, \widetilde{\pi}_1} (s) 
L_{\pi_1, \widetilde{\pi}_2} (s) 
L_{\pi_2, \widetilde{\pi}_1} (s) 
L_{\pi_2, \widetilde{\pi}_2} (s)
\nonumber
\\
&=& 
\Big(\sum_{n=1}^{\infty} 
\frac{|\lambda_{\pi_1} (n^2)|^2}{n^s} \Big)
\Big(\sum_{n=1}^{\infty} 
\frac{\lambda_{\pi_1}(n^2)\bar{\lambda}_{\pi_2}(n^2)}{n^s}\Big)
\Big(\sum_{n=1}^{\infty} 
\frac{\lambda_{\pi_2}(n^2)\bar{\lambda}_{\pi_1}(n^2)}{n^s}\Big)
\Big(\sum_{n=1}^{\infty} 
\frac{|\lambda_{\pi_2} (n^2)|^2}{n^s} \Big)
\nonumber
\\
&=& 
\sum_{n=1}^{\infty} 
\frac{\lambda_{\wedge^2}(n)}{n^s}
\end{eqnarray}
for Re $s$ sufficiently large. Here
$$ 
\lambda_{\wedge^2}(n)
= 
\sum_{\substack{a, b, c, d \geq 1 \\ abcd =n } } 
|\lambda_{\pi_1} (a^2)|^2 
\lambda_{\pi_1} (b^2) 
\bar{\lambda}_{\pi_2} (b^2)
\lambda_{\pi_2} (c^2) 
\bar{\lambda}_{\pi_1} (c^2)  
|\lambda_{\pi_2} (d^2)|^2.
$$

By Lemma 3.1, the four $L$-functions in the product of 
$L_{\wedge^2}(s)$ in \eqref{4L} can be factored as 
\begin{equation}\label{Lij}
L_{\pi_i,\pi_j}(s)
=
L(s,Sym^2\pi_i\times Sym^2\widetilde\pi_j)
U_{i,j}(s)
\end{equation}
for $i,j=1,2$, where $U_{i,j}(s)$ as an Euler product 
converges absolutely for Re $s>\frac12$. Multiplying 
\eqref{Lij} according to \eqref{4L} we can express 
$L_{\wedge^2}(s)$ as the 
product of \eqref{exteriorL} and 
$\prod_{i,j=1}^2U_{i,j}(s)$. Then using Landau's Lemma 2.1 
with $k=2$ and $2A=36$ and Dirichlet convolution we prove 
\eqref{Thm16}.\hfill$\square$

{\it Proof of Theorem 1.7.} From \eqref{notpieta} we can show that $\pi$ and 
$\pi \otimes \eta_{F/\mathbb{Q}}$ are not twisted equivalent 
by any $|\det|^{it}$, $t \in \mathbb{R}$. In fact, if
\begin{eqnarray*}
 \pi \cong \pi \otimes \eta_{F/\mathbb{Q}} \otimes |\det|^{it},
\end{eqnarray*}
for some $ t \in \mathbb{R}^{\times}$, then
\begin{eqnarray*}
 \pi \cong \pi \otimes  |\det|^{i \ell t},
\end{eqnarray*}
which is impossible. Moreover,
\begin{eqnarray}\label{nottwist}
 \pi \otimes \eta_{F/\mathbb{Q}}^{j} 
\ncong \pi \otimes \eta_{F/\mathbb{Q}}^{k} 
\otimes |\det|^{it}
\end{eqnarray}
for any $t \in \mathbb{R}$ when $0 \leq j < k <\ell$. 

By \eqref{LPi} we have the Rankin-Selberg $L$-function
\begin{equation}\label{RSPi}
L(s, \Pi \times \widetilde{\Pi}) 
=
\prod_{j=0}^{\ell-1}   
\prod_{k=0}^{\ell-1} 
L(s, (\pi \otimes \eta_{F/\mathbb{Q}}^{j}) 
\times
(\widetilde\pi \otimes \bar\eta_{F/\mathbb{Q}}^{k})).
\end{equation}
By \eqref{nottwist} we can see on the right hand side of 
\eqref{RSPi}, each diagonal $L$-function ($j=k$) has a simple 
pole at $s=1$ and each off-diagonal $L$-function ($j\neq k$) 
is entire. Consequently, $L(s, \Pi \times \widetilde{\Pi})$ 
has a pole of order $\ell$ at $s=1$ and is analytic elsewhere. 
It also has a standard functional equation. Its complete 
$L$-function is bounded in vertical strips and is holomorphic except for poles at $s=0, 1$ (\cite{JcqSha1}, 
\cite{JcqSha2}, \cite{MglWld}).

The $L$-function in \eqref{RSPi} can be simplified to
\begin{eqnarray}\label{simplified}
L(s, \Pi \times \widetilde{\Pi}) 
=
\prod_{j=0}^{\ell-1}   
L(s, \pi  \times \widetilde{\pi} \otimes 
\eta_{F/\mathbb{Q}}^{j})^\ell.
\end{eqnarray}
Note this $L$-function has nonnegative coefficients
 \begin{eqnarray*}
L(s, \Pi \times \widetilde{\Pi})
=\sum_{n=1}^{\infty}
\frac{\lambda_{\Pi \times \widetilde{\Pi}} (n)}{n^s}, 
\quad 
\lambda_{\Pi \times \widetilde{\Pi}} (n) \geq 0.
\end{eqnarray*}
It converges absolutely for Re $s>1$. Its rank is $4 \ell^2$ 
and hence the sum of coefficients of $s$ in $\Gamma$ functions in its Archimedean factor is
$2 \ell^2$. By Landau's Lemma 2.1,
\begin{eqnarray}\label{asymBC}
\sum_{n \leq x} 
\lambda_{\Pi \times \widetilde{\Pi}} (n) 
= 
x P_{\ell-1} (\log x) 
+O\Big(x^{\frac{4 \ell^2 -1}{4 \ell^2 +1}} 
\log^{\ell-1} x \Big),
\end{eqnarray}
where $P_{\ell-1} $ is a polynomial of degree $\ell-1$.

By \eqref{simplified}, we may express the Fourier coefficients 
$\lambda_{\Pi \times \widetilde{\Pi}} (n)$ in terms of 
$\lambda_{\pi} (n)$. By Goldfeld and Hundley \cite{GldHnd}, 
$\eta_{F/\mathbb{Q}}$ is an idelic lift of a nontrivial Dirichlet character $\chi\!\! \mod \ell$. Therefore by 
\eqref{simplified} 
\begin{eqnarray}\label{lambdaBC}
\lambda_{\Pi \times \widetilde{\Pi}} (n) 
= 
\sum_{\substack{a_{ij} \geq 1 \\ 
\prod_{i=0}^{\ell-1} \prod_{j=0}^{\ell-1} a_{ij} =n}}
\prod_{i=0}^{\ell-1} \prod_{j=0}^{\ell-1} 
|\lambda_{\pi} (a_{ij})|^2 \chi^j (a_{ij}).
\end{eqnarray}
Then \eqref{Thm17} follows from \eqref{asymBC} and 
\eqref{lambdaBC}.
\hfill$\square$

\section{Proof of Theorem 1.8}
\setcounter{equation}{0}

Let $f$ be a self-dual holomorphic cusp form for a 
congruence subgroup of 
$SL_2(\mathbb Z)$. To improve the error term bound in 
\eqref{n2asym}, we factor \eqref{LmUm} further for $m=2$:
\begin{equation}\label{L2zetaU2}
L_2(s)
=
\zeta(s)L(s,Sym^2f)L(s,Sym^4f)U_2(s)
\end{equation}
using the local parameters 
$$
\alpha_p^4,\ \alpha_p^2\times2,\ 1\times3,\ \beta_p^2\times2,
\ \beta_p^4
$$
of the self-dual $L(s,Sym^2f\times Sym^2f)$. 
As in the proof of Theorem 1.2 in \cite{LaoSnk2}, we use 
Perron's formula and shift the contour from 
$\int_{1+\varepsilon-iT}^{1+\varepsilon+iT}$ to 
$\int_{1/2+\varepsilon-iT}^{1/2+\varepsilon+iT}$ to get 
\begin{eqnarray}\label{I1I2I3}
\sum_{n \leq x}
\lambda_f^2(n^2)
&=&
\frac1{2\pi i}
\Big\{
\int_{\frac12 +\varepsilon-iT}^{\frac12 +\varepsilon +\varepsilon+iT}
+
\int_{\frac12 +\varepsilon+iT}^{1+\varepsilon+iT}
+
\int_{1+\varepsilon-iT}^{\frac12 +\varepsilon+\varepsilon-iT}
\Big\}
L_2(s)\frac{x^s}{s}ds 
\nonumber
\\
&&+
{\mathop{Res}\limits_{s=1}}\,
\Big( 
\frac {L_2(s)x^s}{s}\Big)
+O\Big(\frac{x^{1+\varepsilon}}{T}\Big)
\nonumber
\\
&=:& 
I_1+I_2+I_3+{c}_{2}x
+O\Big(\frac{x^{1+\varepsilon}}{T}\Big).
\end{eqnarray}

Since $U_2(s)$ converges absolutely for Re $s>1/2$,
\begin{eqnarray*}
I_1
&\ll& 
x^{\frac12 +\varepsilon}
+
x^{\frac12 +\varepsilon}
\int_1^{T}
\Big|
L\Big(\frac12 +\varepsilon+it, Sym^2f \times Sym^2f \Big)
{U}_2\Big(\frac12 +\varepsilon+it\Big)
\Big| 
\frac{dt}t
\nonumber
\\
&\ll& 
x^{\frac12 +\varepsilon}
+
x^{\frac12 +\varepsilon}
\int_1^{T}
\Big|
L\Big(\frac12 +\varepsilon+it, Sym^2f \times Sym^2f \Big)
\Big| 
\frac{dt}t.
\end{eqnarray*}
By a dyadic subdivision we have
$$
I_1 
\ll 
x^{\frac12 +\varepsilon}
+
x^{\frac12 +\varepsilon}
(\log T) 
\max_{T_1 \leq T}
\Big\{
\frac{1}{T_1}
\int_{\frac{T_1}{2}}^{T_1}
\Big|
L\Big(\frac12 +\varepsilon+it,Sym^2f \times Sym^2f\Big)
\Big|dt \Big\}.
$$
By \eqref{L2zetaU2} and Cauchy's inequality we have
\begin{eqnarray}\label{I1ll}
I_1
&\ll& 
x^{\frac12 +\varepsilon} 
+   
x^{\frac12 +\varepsilon}
(\log T) 
\max_{T_1 \leq T}
\frac{1}{T_1}
\Big(
\max_{\frac{T_1}{2}  \leq t \leq T_1}
\Big|
\zeta\Big(\frac12 +\varepsilon+it\Big)
\Big|\Big) 
\nonumber
\\
&&\times
\Big(
\int_{\frac{T_1}{2}}^{T_1}
\Big|
L\Big(\frac12 +\varepsilon+it,Sym^2f \Big)
\Big|^{2}dt\Big)^{\frac{1}{2}}
\Big(
\int_{\frac{T_1}{2}}^{T_1}
\Big|
L\Big(\frac12 +\varepsilon+it,Sym^4f\Big)
\Big|^{2}dt\Big)^{\frac{1}{2}}.
\end{eqnarray}
Now we need Lemma 2.6 in \cite{LaoSnk2} 
\begin{equation}\label{symbound}
\int_{\frac{T}{2}}^{T}
\Big|
L\Big(\frac12 +\varepsilon+it,Sym^mf \Big)
\Big|^{2}dt
\ll_{f,\varepsilon}
T^{\frac{m+1}2+\varepsilon}
\end{equation}
and a subconvexity bound for $\zeta(s)$ by Bourgain 
\cite{Brg} 
\begin{equation}\label{zetabound}
\zeta\Big(\frac12+\varepsilon+it\Big)
\ll
(1+|t|)^{\frac{53}{342}+\varepsilon}.
\end{equation}
Applying \eqref{symbound} and \eqref{zetabound} to 
\eqref{I1ll} we have
\begin{equation}\label{I1bound}
I_1 
\ll 
x^{\frac12 +\varepsilon}
+
x^{\frac12 +\varepsilon}  
T^{-1+ \frac{53}{342}+ \frac34 +\frac54 +\varepsilon} 
\ll 
x^{\frac12 +\varepsilon}
T^{\frac{395}{342}+\varepsilon}.
\end{equation}

For the integrals $I_2$ and $I_3$ over the horizontal segments, we use \eqref{zetabound} and the convexity bounds for 
$L(s,Sym^2f)$ and $L(s,Sym^4f)$ to get
\begin{eqnarray}\label{I2I3bound}
I_2+I_3
&\ll&
\max_{\frac12 +\varepsilon \leq \sigma \leq 1+\varepsilon}
x^{\sigma}
T^{(\frac{53}{171}+\frac{3}{2}+\frac{5}{2})(1-\sigma)-1+\varepsilon} 
\nonumber
\\
&\ll& 
\Big(\frac{x}{T^{\frac{737}{171}}}\Big)^{1+\varepsilon}T^{\frac{737}{171}-1+\varepsilon} + \Big(\frac{x}{T^{\frac{737}{171}}}\Big)^{\frac12 +\varepsilon}T^{\frac{737}{171}-1+\varepsilon}
\nonumber
\\
&\ll& 
\frac{x^{1+\varepsilon}}{T}
+
x^{\frac12 +\varepsilon }
T^{\frac{395}{342}+\varepsilon}.
\end{eqnarray}
By \eqref{I1I2I3}, \eqref{I1bound} and \eqref{I2I3bound} we get
\begin{equation}\label{xTbound}
\sum_{n \leq x}
\lambda_f^2(n^2)
=
c_{2}x+
O\Big(\frac{x^{1+\varepsilon}}{T}\Big)
+O\Big(x^{\frac12 +\varepsilon}T^{\frac{395}{342}+\varepsilon}\Big).
\end{equation}
Taking $T=x^{\frac{171}{737}}$ in \eqref{xTbound}, we finally 
get \eqref{n2asymholo}.\hfill$\square$

{\bf Acknowledgments.} This work was completed when the first 
author visited the University of Iowa, supported by the 
International Coorporation Program sponsored by the Shandong 
Provincial Education Department.

\end{document}